\documentclass{crmp-l}

\usepackage{graphics}

\newtheorem{theorem}{Theorem}[section]

\theoremstyle{definition}
\newtheorem{definition}[theorem]{Definition}

\theoremstyle{remark}

\numberwithin{equation}{section}

\newcommand{\abs}[1]{\lvert#1\rvert}

\newcommand{\FF}{\mathbb{F}}
\newcommand{\PP}{\mathbb{P}}
\newcommand{\M}{\mathbb{M}}
\newcommand{\Z}{\mathbb{Z}}
\newcommand{\inc}[1]{\ensuremath{\mathit{Inc}(#1)}}
\newcommand{\udl}[1]{\underline{#1}}

\begin{document}

\title{Deflating infinite {C}oxeter groups to finite groups}

\author{Christopher S. Simons}
\address{Department of Mathematics, Rowan
University, Glassboro, New Jersey, 08028, U.S.A.}
\email{simons@alumni.princeton.edu}
\thanks{The author was supported in part by an NSERC postdoctoral
fellowship.}

\date{\today}

\maketitle

\section{Introduction} 

In \cite{implies26} it is shown, using the Ivanov-Norton theorem
\cite{iva165,nor165}, that the Bimonster, $\M\wr 2$, is presented by the
Coxeter relations of the 26 node incidence graph of the projective plane of
order 3 along with the additional relations that all free 12-gons of this
diagram generate symmetric groups $S_{12}$.  The purpose of this paper is
to consider an easier version of this presentation by replacing the 26 node
incidence graph of the projective plane of order 3, $\inc{\PP_3}$, by the
14 node incidence graph of the projective plane of order 2, $\inc{\PP_2}$. 
The maximum free $n$-gons of $\inc{\PP_2}$ are 8-gons so the additional
non-Coxeter relations are that all free 8-gons generate symmetric groups
$S_8$.  We prove that the resulting group is $O^{-}_8(2){:}2$ (in the Atlas
\cite{atlas} notation).

This new presentation is in some sense more satisfying than the related
Monstrous presentation as it is possible to understand it in an elementary
and self contained fashion.  We hope to use it and other similar
presentations to explain the finite simple Monster group in a nonsporadic
manner.  That $O^{-}_8(2){:}2$ acts as a pseudomonster is not entirely new
as it is the unique nonsymmetric $S(3,8)$ group while the $\M\wr 2$ is the
unique nonsymmetric $S(5,12)$ group \cite{con93}.

To avoid unnecessary complications, in Section \ref{results} we restrict
our attention to the $O^{-}_8(2){:}2$ case.  Comments on the Bimonster and
general cases are left to Section \ref{epilogue}. A useful basic reference
on Coxeter and reflection groups is \cite{rgacg}.  We also provide
some relevant references for the motivating Monstrous case.

\section{Results} 
\label{results}

\begin{figure}
\centering\includegraphics{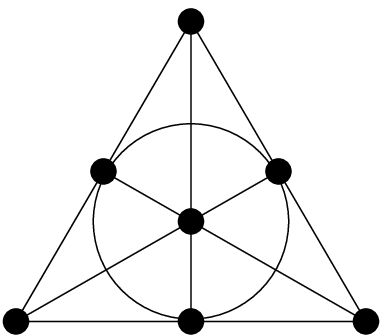}
\caption{$\PP_2$}
\label{proj2}
\end{figure}

\begin{figure}
\centering\includegraphics{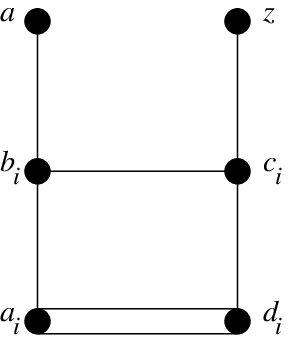}
\caption{$\inc{\PP_2}$}
\label{incproj2}
\end{figure}

The projective plane of order 2, $\PP_2$, consists of 7 points and 7 lines
as shown in Figure \ref{proj2}. Each line consist of 3 points, any two
lines uniquely determine a point and there is a duality between points and
lines. The incidence graph of the projective plane of order 2,
$\inc{\PP_2}$, therefore has 14 nodes (one for each of the 7 points and 7
lines of $\PP_2$).  Two nodes are joined exactly when one is a point, one
is a line and the point lies on the line.  The graph $\inc{\PP_2}$ has
valence 3 and is shown in Figure \ref{incproj2}.  The indices $i$ range
over $\{ 1, 2, 3\}$ so that some of the vertices of the figure correspond
to 3 nodes of $\inc{\PP_2}$.  Single lines in the figure indicate that two
of these nodes are joined just if they share the same index.  Double lines
indicate that two of these nodes are joined just if their indices differ.

To $\inc{\PP_2}$ or any graph $\Gamma$ we can associate a Coxeter group
$H$ generated by the nodes and subject to the graph's Coxeter relations. 
The Coxeter relations are that all nodes are involutions ($\alpha^2 =
1$), that any two unjoined nodes commute ($(\alpha\beta)^2 = 1)$)  and
that the product of any two joined nodes has order 3 ($(\alpha\beta)^3 =
1)$). 

It is important for us that all Coxeter groups are isomorphic to reflection
groups.  This allows us to work geometrically.  Each generating involution
corresponds to the reflection, $x \longmapsto x - (x,\alpha) \alpha$, in a
root $\alpha$ (we use the same name for both the root and the associated
reflection).  The inner products of these roots are such that
$(\alpha,\alpha) = 2$ for all roots $\alpha$, $(\alpha,\beta) = 0$ for
$\alpha, \beta$ unjoined and $(\alpha,\beta) = -1$ for $\alpha, \beta$
joined. 

The $\inc{\PP_2}$ Coxeter group is infinite.  In order to present the
finite group $O^{-}_8(2){:}2$ we must adjoin some additional relations.
We make use of the following convention introduced in \cite{implies26}.

\begin{definition} 
Let $\Z^m{:}G$ be a group (usually an affine Coxeter group).  To deflate
this group is to impose the relations that make the translations $\Z^m$
trivial.  Similarly to biflate, triflate or $k$-flate $\Z^m{:}G$ is to
make the translations have order 2, 3 or $k$ respectively with the result
that $\Z^m$ becomes $2^m$, $3^m$ or $k^m$.  Usually we view $\Z^m{:}G$ as
a subgroup of a group $H$.  The new relations are then imposed on H. 
\end{definition}

$\Z^m$ is the free Abelian group of rank $m$.  The group notation $A{:}B$
is used for a semidirect product of $A$ and $B$. Note that $A$ will be a
normal subgroup of $A{:}B$.  Below we use lower case letters for spherical
Dynkin diagrams and upper case letters to the associated affine Dynkin
diagrams.

We deflate affine $A_7$ Coxeter groups $\Z^7{:}S_8$.  Since the affine
$A_7$ diagram is just a free $8$-gon (a graph with 8 nodes joined precisely
as the vertices of an octagon), we say that we are deflating a free
8-gon.  For example, if $\alpha_0, \ldots ,\alpha_7$ are the involutions
of an 8-gon then we can deflate this 8-gon by adding the relation
$\alpha_0 = \alpha_1^{\alpha_2\cdots\alpha_{7}}$.  This should be clear in
a moment. 

We now state the main result of this paper.

\begin{theorem}  
\label{theorem1}
If we deflate all free 8-gons of the $\inc{\PP_2}$ Coxeter group we
obtain the finite group $O^{-}_8(2){:}2$.
\end{theorem}

In general it is very difficult to identify a group from its presentation.
To prove Theorem \ref{theorem1} we therefore take an indirect approach.

\begin{figure}
\centering\includegraphics{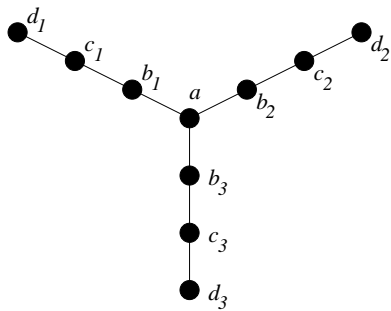}
\caption{$Y_{333}$ diagram}
\label{y333diagram}
\end{figure}

We start with a (small) graph $\Gamma_0$: the $Y_{333}$ digram (Figure
\ref{y333diagram}).  We consider a group $G$ satisfying its Coxeter
relations.  We then extend $\Gamma_0$ by repeatedly adjoining affine $A_7$
extending nodes.  The closure of $\Gamma_0$ under this extension is denoted
$\Gamma$.  Of course $G$ satisfies the Coxeter relations of $\Gamma$. 

Adjoining $A_7$ extending nodes is straightforward.  Let the nodes
$\alpha_1, \ldots, \alpha_7$ form an $a_7$ subdiagram (7 nodes forming a
chain). Its Coxeter group is $S_8$ and the nodes can be viewed as the
transpositions $(0,1), \dots, (6,7)$ of $S_8$.  The extending node
$\alpha_0$ is then the transposition $(7,0)$ and can be written as
$\alpha_1^{\alpha_2\cdots\alpha_{7}}$.  In order to make use of $\alpha_0$
in further $A_7$ extensions we must determine its Coxeter relations with
other nodes of $\Gamma$.  We do this by treating $G$ as (the quotient of)
a reflection group.  We carefully choose root vectors $\alpha_1, \ldots,
\alpha_7$.  These must have norm 2, inner product 0 when unjoined and
inner product $-1$ when joined.  (Ensuring inner product $-1$ requires
some care since while the sign of a root does not effect the corresponding
reflection element it is important here.  The choice of signs does not
have to be consistent between different $a_7$ subdiagrams.)  We then have
by Coxeter theory that the extending root $\alpha_0$ is $-(\alpha_1 +
\cdots + \alpha_7)$.  Inner products can determine the Coxeter relations
between the new node $\alpha_0$ and other nodes $\alpha$ of $\Gamma$.  If
$(\alpha_0,\alpha) = 0$ then $\alpha_0$ and $\alpha$ are unjoined.  If
$(\alpha_0,\alpha) = \pm 1$ then $\alpha_0$ and $\alpha$ are joined. 

Two roots corresponding to the same element of $G$ are said to be equivalent
and they share the same node of $\Gamma$.

We prove that if $\Gamma$ (for $\Gamma_0$ still the $Y_{333}$ diagram) has
no more than 14 distinct nodes then $G$ is $O^{-}_8(2){:}2$.  Upon
examination of the relations used in this proof we then get Theorem
\ref{theorem1}. 

\begin{theorem}[14 implies $O^{-}_8(2){:}2$]
\label{theorem2}

Let $G$ be a group (of order greater than 2) generated by and satisfying
the Coxeter relations of the $Y_{333}$ diagram.  If after closure of the
generating $Y_{333}$ diagram by adjoining $A_7$ extending nodes there are
no more than 14 nodes, then $G$ is $O^{-}_8(2){:}2$. 

\end{theorem} 

Insisting that there are no more than 14 nodes imposes some non-Coxeter
relations, the equivalence of certain roots, on the $Y_{333}$ infinite
Coxeter group with the result that we get a presentation for
$O^{-}_8(2){:}2$.  In fact we show that $\Gamma$ has exactly 14 nodes and
is $\inc{\PP_2}$.  We do assume for now that $G$ is not trivial (of order
1 or 2). This is justified since by its construction the group
$O^{-}_8(2){:}2$ will satisfy the conditions.

We start with $\Gamma_0$ being the $Y_{333}$ diagram and close under $A_7$
extension.  In order to proceed we use the following coordinate system. 
Often we will use $+$ to denote $1$ and $-$ to denote $-1$. 

We have a space of 13 coordinates
\begin{equation}
\begin{array}{rrrrr}
  a & b & c & d &       \\
  e & f & g & h & t \\
  i & j & k & l &
\end{array}
\text{ with quadratic form } a^2 + \dotsb + l^2 - t^2.
\label{syseqn}
\end{equation}

\begin{figure}
\centering\includegraphics{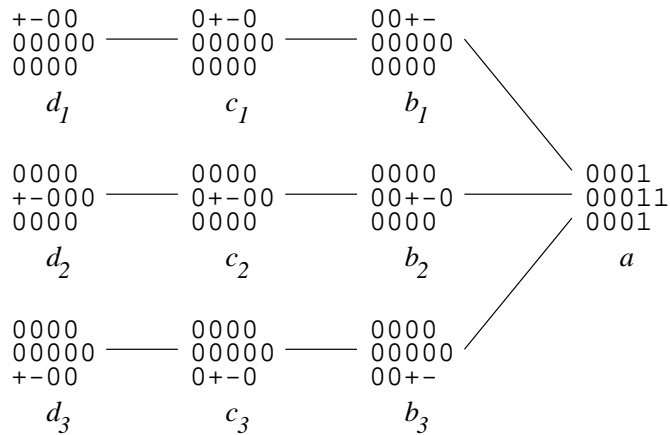}
\caption{The fundamental roots of $Y_{333}$}
\label{sys}
\end{figure}

In this system the 10 original roots of $Y_{333}$ are as indicated in
Figure \ref{sys}. 

All the vectors satisfy the following relations:
\begin{align}
a+b+c+d=t \nonumber \\
e+f+g+h=t  \\
i+j+k+l=t . \nonumber
\end{align}
Thus $t$ is redundant, so can be omitted.

We now extend the $Y_{333}$ diagram by adjoining $A_7$ extending roots to
obtain $\Gamma$.  We do this by finding free $a_7$ subdiagrams and adding
the $A_7$ extending nodes.  We stress that in all cases used the $a_7$
roots satisfy the standard inner product conditions.  The choice of signs
of the $a_7$ roots is important.  If $\alpha$ is a root then we use
$\udl{\alpha}$ to denote the negative of $\alpha$.  The signed sum of the
$a_7$ roots is the $A_7$ extending root. 

The $Y_{333}$ diagram has 10 nodes.

We use the subdiagram $d_1$---$c_1$---$b_1$---$a$---$b_2$---$c_2$---$d_2$
to get
\begin{equation}
a_3 =
\begin{array}{rrrrr}
  1 & 0 & 0 & 0 &    \\
  1 & 0 & 0 & 0 &  1 \\
  0 & 0 & 0 & 1 & 
\end{array}. 
\end{equation}

Checking inner products we see that $a_3$ is joined to $b_3$ in addition to
$d_1$ and $d_2$.  It is unjoined with the other nodes of $\Gamma$.  

By the obvious $S_3$ symmetry of the the $Y_{333}$ diagram we similarly
obtain $a_1$ and $a_2$.  

\begin{figure}
\centering\includegraphics{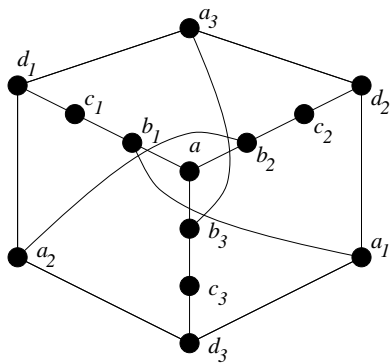}
\caption{$a_1$, $a_2$, $a_3$}
\label{a123}
\end{figure} 

So far we have 13 roots as shown in Figure \ref{a123}.  By checking inner
products with the original 10 roots of $Y_{333}$ we find that the roots
are distinct.  So under the conditions of Theorem \ref{theorem2} there is
only one more node of $\Gamma$ still to be found. 

We now consider the subdiagram
$\udl{c_1}$---$\udl{d_1}$---$a_2$---$\udl{d_3}$---$a_1$---$\udl{d_2}$---$\udl{c_2}$.
The extending root is
\begin{equation}
z_3 =
\begin{array}{rrrrr}
  0 & 0 & 1 & 1 &    \\
  0 & 0 & 1 & 1 &  2 \\
  1 & 1 & 0 & 0 &
\end{array}.
\end{equation}

Checking inner products we see that $z_3$ is joined to $c_1$, $c_2$, $c_3$
and is unjoined with the other nodes of $\Gamma$.

We can similarly obtain $z_1$ and $z_2$.  These have the same Coxeter
relations as $z_3$.  So by the conditions of Theorem \ref{theorem2} we
have that the $z_i$ are equivalent, $z_1 \equiv z_2 \equiv z_3$:
\begin{equation}
\begin{array}{rrrrr}
  1 & 1 & 0 & 0 &    \\
  0 & 0 & 1 & 1 &  2 \\
  0 & 0 & 1 & 1 &
\end{array}
\equiv
\begin{array}{rrrrr}
  0 & 0 & 1 & 1 &    \\
  1 & 1 & 0 & 0 &  2 \\
  0 & 0 & 1 & 1 &
\end{array}
\equiv
\begin{array}{rrrrr}
  0 & 0 & 1 & 1 &    \\
  0 & 0 & 1 & 1 &  2 \\
  1 & 1 & 0 & 0 &
\end{array}.
\label{extrarels}
\end{equation}
We call this common node $z$.

We now have 14 distinct nodes and it is easily checked that the graph
obtained is $\inc{\PP_2}$.  It is not possible to make any further $A_7$
extensions.  Therefore $\Gamma = \inc{\PP_2}$.

We now identify the group $G$.  We do so by explicit enumeration of the
root elements of $G$.  As this has been described in detail in
\cite{con93} we merely sketch the proof here here.  We use the $Y_{333}$
reflection group along the relation of equation \eqref{extrarels} and find
that the 136 root elements represented by 

\begin{equation}
\begin{array}{rrccr} 
0&0&+&-& \\ 
0&0&0&0&0 \\ 
0&0&0&0& \\
\multicolumn{5}{c}{(18)} 
\end{array}\quad 
\begin{array}{rrrrr} 
0&0&0&1& \\
0&0&0&1&1 \\ 
0&0&0&1& \\ 
\multicolumn{5}{c}{(64)} 
\end{array}\quad
\begin{array}{rrrrr} 
0&0&1&1& \\ 
0&0&1&1&2 \\ 
0&0&1&1& \\
\multicolumn{5}{c}{(54)} 
\end{array} 
\end{equation} 
form a conjugacy class of $G$.  [The elements are obtained from the roots
shown by permuting all four coordinates of each row and permuting the
three rows.  Parentheses show the number of distinct group elements
obtained.]

These roots are then transformed by $G$ in exactly the same way that
$O^{-}_8(2){:}2$ permutes its root elements.  [To see this view the 12
coordinates $a, \ldots, l$ modulo 2.  The roots span a 10 dimensional
subspace of $\FF_2^{12}$.  After applying the relations of equation
\eqref{extrarels} we get an 8 dimensional space $\FF_2^8$.  The
bilinear form inherited from the quadratic form (equation \eqref{syseqn})
has Witt defect 1.]   It follows the resulting group is a central
extension of $O^{-}_8(2){:}2$.  The multiplier of the simple
group $O^{-}_8(2)$ is known to be 1 \cite{atlas}, therefore the extension
is trivial and $G$ is $O^{-}_8(2){:}2$.

This proves Theorem \ref{theorem2}.  We note again that Theorem
\ref{theorem1} quickly follows since all the relations used in this proof
are implied by the conditions of Theorem \ref{theorem1} and all
relations of Theorem \ref{theorem1} hold in $O^{-}_8(2){:}2$. 

\section{Epilogue}
\label{epilogue}

For reference we state the main theorems of \cite{implies26}.

\begin{theorem}
If we deflate all free 12-gons of the $\inc{\PP_3}$ Coxeter group we
obtain the Bimonster $\M\wr 2$.
\end{theorem}

\begin{theorem}[26 implies the Bimonster] Let $G$ be a group (of order
greater than 2) generated by and satisfying the Coxeter relations of the
$\M_{666}$ diagram.  If after closure of the generating $\M_{666}$ diagram
by adjoining $A_{11}$ extending nodes there are no more than 26 nodes, then
$G$ is the Bimonster $\M\wr 2$.  \end{theorem}

$\inc{\PP_3}$ is the (26 node) incidence graph of the projective plane of
order 3.  The $\M_{666}$ diagram is the (16 node) $Y_{555}$ diagram. 

The proofs of these theorems are similar those given in this paper.  The
principal complication is that the root enumeration would involve
$\abs{\M}$ (almost $10^{54}$) root elements and therefore much less
elementary methods \cite{iva165, nor165, iva99} are required to identify
the group. 

We also include a table (Table \ref{game1}) from \cite{implies26}.
In this table we list the groups obtained by deflating all free
$n$-gons of the Coxeter groups of certain graphs.

\begin{table}
\centering{\begin{tabular}{|r|r|c|c|}
\hline
$n$-gons & \# of nodes & graph & group \\
\hline
12 & 26 & $ \inc{\PP_3} $ & $ \M\wr 2 $ \\
8  & 14 & $ \inc{\PP_2} $ & $ O^{-}_8(2){:}2 $ \\
6  & 10 & Petersen & $ O^{-}_6(2){:}2 \cong O_5(3){:}2 $ \\
6  &  8 & cube = $\inc{\text{tetrahedron}}$ & $ O_5(3)\times 2 \cong
O^{-}_6(2) \times 2 $ \\
\hline
\end{tabular}}
\caption{deflating $n$-gons in graphs to get groups}
\label{game1}
\end{table}

We are very interested in which other finite groups can be obtained as
deflations of Coxeter groups.

\bibliographystyle{amsalpha}

\begin{thebibliography}{10}

\bibitem [1]{con93} J.~H. Conway.  
\newblock From hyperbolic reflections to finite groups.  
\newblock In L.~Finkelstein and W.~M. Kantor, editors, {\em
Groups and Computation}, number 11 in {DIMACS} Series in Discrete
Mathematics and Theoretical Computer Science, pages 41--51, American
Mathematical Society, 1993. 

\bibitem [2]{atlas} J.~H. Conway, R.~T. Curtis, S.~P. Norton, R.A. Parker,
and R.~A. Wilson. 
\newblock {\em Atlas of Finite Groups}. 
\newblock Oxford University Press, 1985.

\bibitem [3]{cns} J.~H. Conway, S.~P. Norton, and L.~H. Soicher.  
\newblock The {B}imonster, the group {$Y_{555}$}, and the projective plane
of order $3$.  
\newblock In M.~C. Tangara, editor, {\em Computers in Algebra}, number 111
in Lecture Notes in Pure and Applied Mathematics, pages 27--50. Marcel 
Dekker, 1988.

\bibitem [4]{conprit}
J.~H. Conway and A.~D. Pritchard.
\newblock Hyperbolic reflections for the {B}imonster and
$3\mathit{Fi}_{24}$.
\newblock In Liebeck and Saxl \cite{lms165}, pages 24--45.

\bibitem [5]{implies26} J.~H. Conway and C.~S. Simons.  
\newblock 26 implies the {M}onster.  
\newblock Journal of Algebra, in press.  

\bibitem [6]{relations} J.~H. Conway and C.~S. Simons.  
\newblock Relations in {$\M_{666}$}. 
\newblock In R.~T. Curtis and R.~A. Wilson, editors, {\em The
Atlas 10 Years On}, number 249 in London Mathematical Society Lecture Note
Series, pages 27--38. Cambridge University Press, 1998. 

\bibitem [7]{rgacg} J.~E. Humphreys.  
\newblock {\em Reflection Groups and {C}oxeter Groups}. 
\newblock Number~29 in Cambridge Studies in Advanced Mathematics. Cambridge 
University Press, 1990.

\bibitem [8]{iva165} A.~A. Ivanov.  
\newblock A geometric characterization of the {M}onster. 
\newblock In Liebeck and Saxl \cite{lms165}, pages 46--62.

\bibitem [9]{iva99} A.~A. Ivanov.  
\newblock {$Y$}-groups via transitive extension. 
\newblock {\em Journal of Algebra}, 218:412--435, 1999.

\bibitem [10]{lms165} M.~Liebeck and J.~Saxl, editors.  
\newblock {\em Groups, Combinatorics and Geometry}, number 165 in London 
Mathematical Society Lecture Note Series. Cambridge University Press,
1992.

\bibitem [11]{nor165} S.~P. Norton.
\newblock Constructing the {M}onster.
\newblock In Liebeck and Saxl \cite{lms165}, pages 63--76.

\bibitem [12]{simons2} C.~S. Simons.  
\newblock Monster roots.  
\newblock In J.~Ferrar and K.~Harada, editors, {\em The Monster and Lie 
Algebras}, number~7 in Ohio State University Mathematical Research 
Institute Publications, pages 127--146. Walter de Gruyter, 1998.

\end{thebibliography}

\end{document}